\documentclass[aip,
amsmath,amssymb,
preprint,%
numerical,
]{revtex4-1}
\usepackage{graphicx}
\usepackage{dcolumn}
\usepackage{bm}
\usepackage{hyperref}
\usepackage[utf8]{inputenc}
\usepackage[T1]{fontenc}
\usepackage{mathptmx}
\usepackage{etoolbox}
\newcommand{\RNum}[1]{\uppercase\expandafter{\romannumeral #1\relax}}
\usepackage{array,xcolor,colortbl}
\usepackage{multirow}
\makeatletter
\def\@email#1#2{%
	\endgroup
	\patchcmd{\titleblock@produce}
	{\frontmatter@RRAPformat}
	{\frontmatter@RRAPformat{\produce@RRAP{*#1\href{mailto:#2}{#2}}}\frontmatter@RRAPformat}
	{}{}
}%
\makeatother

\begin{document}
	
	\title{Effect of oblique irradiation on the onset of thermal phototactic bioconvection in non-scattering medium}
	
	\author{S. K. Rajput}
	\altaffiliation[Corresponding author: E-mail: ]{shubh.iiitj@gmail.com.}
	\author{M. K. Panda}%
	\affiliation{$^1$ Department of Mathematics, PDPM Indian Institute of Information Technology Design and Manufacturing, Jabalpur 482005, India.
	}%
	

	\begin{abstract}
		
		The linear stability of a suspension of phototactic algae is investigated numerically with particular emphasis on the effects of the angle of
		incidence of the illuminating oblique collimated irradiation with thermal effects. The suspension is illuminated by the oblique collimated irradiation from the top and heated/cooled from the bottom. The linear stability analysis shows that the suspension becomes more unstable as the angle of incidence increases.
		
	\end{abstract}
	
	
	\maketitle

	
	\section{INTRODUCTION}
	
	The phenomena of spontaneous pattern formation in suspensions of randomly, but on average upwardly swimming microorganisms are known as bioconvection~\cite{20platt1961,21pedley1992,22hill2005,23bees2020,24javadi2020}. The swimming microorganisms participating in bioconvection are up to $10\%$ denser than the medium (water here), and they are mostly algae and bacteria. Usually, bioconvection patterns are observed in the laboratory in shallow suspensions. However, these have also been found in situ in micro patches of zoo plankton.  Also, bioconvection patterns disappear when the microorganisms stop swimming. Nevertheless, examples of pattern formation are also found where up-swimming and higher density are not involved.  Microorganisms can adjust their swimming in response to several different environmental stimuli, which can cause them to aggregate in speciﬁc regions. These responses are known as taxes. Paradigmatic examples of taxes include gravitaxis, chemotaxis, gyrotaxis, and phototaxis. Gravitaxis refers to the response of microorganisms to gravity or acceleration, and negative gravitaxis is the swimming opposite of gravity. Responses to chemical gradients can lead to chemotaxis. Gyrotaxis is termed as the directed swimming of a bottom-heavy microorganism due to the balance between a torque due to gravity and viscous torque arising from local shear ﬂows. Positive (negative) phototaxis is termed as the directed movement toward (away from) the dim (bright) light source, while start/stop swimming behaviour is observed in the photophobic response. The process by which the algae absorb the light incident on them directly from above and produce shadows below them is called self-shading (or shading). This manuscript is restricted to phototaxis and self-shading only.	
	
	Experimental observations have revealed that the bioconvection patterns are modiﬁed by illuminating sources of different types, for example, oblique and/or vertical collimated irradiation. The steady patterns in suspensions of (phototactic) algae may be destroyed due to strong light or their formation may be prevented in well-stirred cultures. The shape, size, symmetry, and/or scale of the bioconvection patterns can also be modiﬁed by illumination. The reason for these changes here is twofold: First, the photosynthetic pigments (e.g., chlorophyll and carotenoid) of the motile algae aid in obtaining energy during their photosynthesis. Since the motile algae are strongly phototactic, they move toward (or away from) the weak (or strong) light source via a light-seeking (or light-avoiding) behavior when $G < G_c$ (or $G > G_c$). Thereby, the cells tend to accumulate at optimal places ($G=G_c$) in their local environment. The absorption of light by the algae may be the second reason for the changes in bioconvection patterns.
	
	The study of bioconvection has garnered significant attention across various fields, particularly focusing on the interaction between thermal and phototactic factors in microorganism suspensions. Kuznetsov~\cite{51kuznetsov2005thermo} delved into bio-thermal convection within suspensions of oxytactic microorganisms, while Alloui et al.\cite{52alloui2006stability} investigated suspensions of mobile gravitactic microorganisms. Nield and Kuznetsov\cite{53nield2006onset} utilized linear stability analysis to explore the onset of bio-thermal convection in suspensions of gyrotactic microorganisms, and Alloui et al.\cite{54alloui2007numerical} scrutinized the impact of bottom heating on the onset of gravitactic bioconvection in a square enclosure. Taheri and Bilgen\cite{55taheri2008thermo} investigated the effects of bottom heating or cooling in a vertically oriented cylinder with stress-free sidewalls. Kuznetsov~\cite{56kuznetsov2011bio} developed a theoretical framework for bio-thermal convection in suspensions containing two species of microorganisms. Saini et al.\cite{57saini2018analysis} explored bio-thermal convection in suspensions of gravitactic microorganisms, while Zhao et al.\cite{57zhao2018linear} utilized linear stability analysis to examine the stability of bioconvection in suspensions of randomly swimming gyrotactic microorganisms heated from below.
	
	In the domain of phototactic bioconvection, Vincent and Hill~\cite{12vincent1996} laid the groundwork with seminal research, examining the impact of collimated irradiation on an absorbing (non-scattering) medium. Building upon this foundation, Ghorai and Hill~\cite{10ghorai2005} extended the inquiry into the behavior of phototactic algal suspensions in two dimensions, without considering scattering effects. Subsequently, Ghorai et al.\cite{7ghorai2010} and Ghorai and Panda\cite{13ghorai2013} delved into the effects of light scattering, both isotropic and anisotropic, under normal collimated irradiation. Panda and Ghorai~\cite{14panda2013} proposed a model for an isotropically scattering medium in two dimensions, yielding results divergent from those of Ghorai and Hill~\cite{10ghorai2005} due to the inclusion of scattering effects. Panda and Singh~\cite{11panda2016} explored phototactic bioconvection in two dimensions, confining a non-scattering suspension between rigid sidewalls. Additionally, Panda et al.\cite{15panda2016} examined the impact of diffuse irradiation, combined with collimated irradiation, in an isotropic scattering medium, while Panda\cite{8panda2020} investigated an anisotropic medium. Recognizing natural environmental conditions where sunlight strikes the Earth's surface at oblique angles, Panda et al.\cite{16panda2022} studied the effects of oblique collimated irradiation on the onset of phototactic bioconvection. In a recent investigation, Panda and Rajput\cite{41rajput2023} explored the impacts of diffuse irradiation along with oblique collimated irradiation on a uniformly scattering suspension. Furthermore, Rajput and Panda~\cite{rajput2024effect} investigated the effect of scattered/diffuse flux on the onset of phototactic bioconvection in the absence of collimated flux. However, no study on the onset of thermal phototactic bioconvection that incorporates the effects of oblique collimated irradiation on an algal suspension has been hitherto carried out. Therefore, the effects of oblique collimated irradiation on thermal bioconvection are investigated in the same vicinity.	
	
	The organization of the manuscript is as follows: The mathematical formulation of the proposed phototaxis model is described ﬁrst.
	The basic steady state is derived next, and then, it is perturbed by
	inﬁnitesimal disturbances. The linear stability problem to the basic
	steady state is derived next and solved numerically by taking two different cases (i.e., stress-free and rigid upper surface). Finally, the results of the numerical study are summarized and their novelty is addressed.	
	
	
	\begin{figure}[!htbp]
		\centering
		\includegraphics[width=12cm]{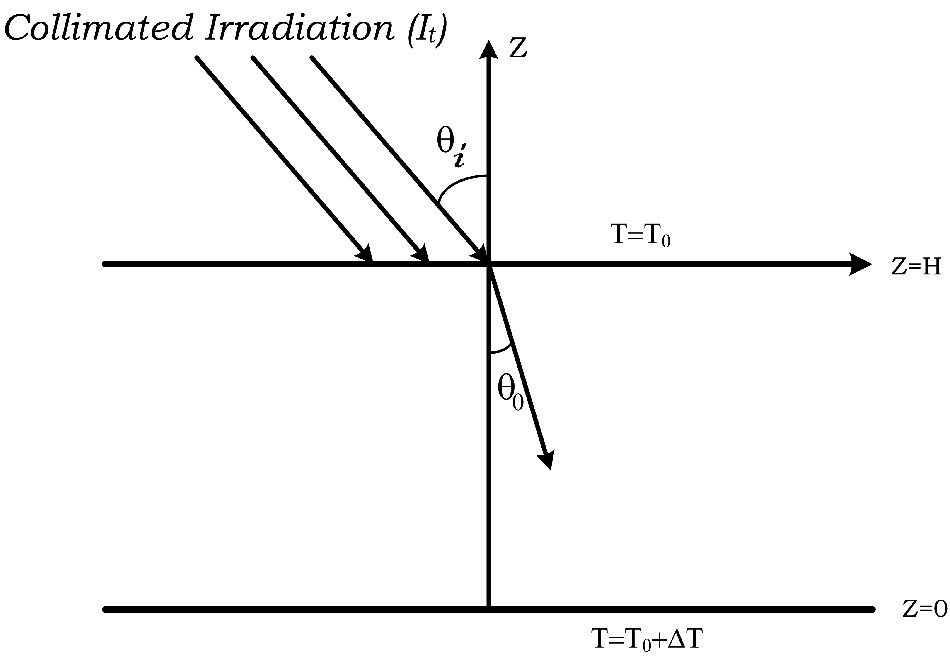}
		\caption{\footnotesize{Oblique collimated irradiation on the upper surface of an absorbing and non-scattering algal suspension with thermal effect on the lower surface.}}
		\label{fig1}
	\end{figure}
	
	\section{MATHEMATICAL FORMULATION}
	
	Consider the motion in a dilute suspension of phototactic algae within a layer of ﬁnite depth $H$. Here, an oblique collimated irradiation illuminates the suspension from above and strikes it at a ﬁxed off-normal angle $\theta_i$. We choose a rectangular Cartesian coordinate system where the yz-plane is the plane of incidence for the oblique collimated irradiation. The angle of refraction $\theta_0$ in which the collimated beam propagates across the water is determined by using Snell’s law, that is, $\sin\theta_i=n_0\sin\theta_0$, Here, $n_0$ denotes the refractive index of water and the estimated value for it is 1.333 approximately. Since the index of refraction of algae is not the same as that of water, the light incident on the algae across the suspension is absorbed by them and scattered thereafter. For simplicity, the effects of scattering by algae have been neglected in the present study. Let $I(\boldsymbol{x},\boldsymbol{s})$ represent the radiation intensity propagating in the (unit) direction $\boldsymbol{s}=\cos\theta\hat{z}+\sin\theta(\cos\phi\hat{x}+\sin\theta\hat{y})$ at a position $\boldsymbol{x}$ across the algal suspension, where $\boldsymbol{x}$ is measured relative to a rectangular cartesian coordinate system with the $z$ axis vertically up. Here, $\theta$ denotes the polar angle (measured from the $z$ axis) and $\phi$ denotes the azimuthal angle (measured in between the projection of the radiation intensity onto xy- plane and the $x$ axis) describing the unit
	vector $$\boldsymbol{x}$$ (in spherical polar coordinate system).
	
	\section{Phototaxis coupled with self-shading in suspensions of algae}
	We assume here that the medium across the algal suspension is absorbing and non-scattering similar to Vincent and Hill. To calculate light intensity proﬁles, the radiative transfer equation (hereafter referred to as RTE) is given by
	
	\begin{equation}\label{1}
	\boldsymbol{s}\cdot\nabla I(\boldsymbol{x},\boldsymbol{s})+\kappa I(\boldsymbol{x},\boldsymbol{s})=0
	\end{equation}
	where $\kappa$ is the absorption coefﬁcient.
	
	The collimated intensity at location $\boldsymbol{x}_b=(x,y,H)$ (i.e., the top boundary surface) is expressed by using a Dirac delta function as
	
	\begin{equation*}
		I(\boldsymbol{x}_b,\boldsymbol{s}) = I_t\delta(\boldsymbol{s}-\boldsymbol{s}_0).
	\end{equation*}
	Here, $I_t$ is the magnitude of oblique collimated irradiation and $\boldsymbol{s}_0=\sin(\pi-\theta_0)\cos\phi\hat{x}+\sin(\pi-\theta_0)\sin\phi\hat{y}+\cos(\pi-\theta_0)\hat{z}$ is the incident direction deﬁned in spherical polar coordinates. Also, $\hat{x}$, $\hat{y}$, and $\hat{z}$ are unit vectors along the $x$, $y$, and $z$ axes. The Dirac-delta function
	$\delta$ satisﬁes
	
	\begin{equation*}
		\int_0^{4\pi}f(\boldsymbol{s})\delta(\boldsymbol{s}-\boldsymbol{s}_0)=f(\boldsymbol{s}_0).
	\end{equation*}
	
	The absorption coefﬁcient is linearly proportional to the concentration $n$ so that $\kappa=\varkappa n$, and thus, the RTE in a non-scattering and absorbing suspension becomes
	
	\begin{equation}\label{2}
	\boldsymbol{s}\cdot\nabla I(\boldsymbol{x},\boldsymbol{s})+\varkappa nI(\boldsymbol{x},\boldsymbol{s})=0
	\end{equation}
	
	The total intensity, $G(\boldsymbol{x})$; at a point $\boldsymbol{x}$ in the medium is
	
	\begin{equation}\label{3}
	G(\boldsymbol{x})=\int_0^{4\pi}I(\boldsymbol{x},\boldsymbol{s})d\Omega.
	\end{equation}
	
	The radiative heat ﬂux, $\boldsymbol{q}(\boldsymbol{x})$, at a point $\boldsymbol{x}$ in the medium is
	
	\begin{equation}\label{4}
	\boldsymbol{q}(\boldsymbol{x})=\int_0^{4\pi}I(\boldsymbol{x},\boldsymbol{s})\boldsymbol{s}d\Omega.
	\end{equation}
	
	Let $\boldsymbol{p}$ be the unit vector pointing in the swimming direction of a cell. The mean swimming direction, $<P>$, is deﬁned as the ensemble average of the swimming direction and p for all the cells in a small volume. The swimming speed does not depend on the illumination, position, time, and direction of many species of microorganisms. It is
	assumed that the cells swim at the same speed relative to the ﬂuid. We denote the ensemble average swimming speed by $W_c$. The average swimming velocity is, thus,
	
	\begin{equation}\label{5}
	<\boldsymbol{p}>=M(G)\hat{\boldsymbol{z}}.
	\end{equation}
	
	Here, $M(G)$ is the phototaxis function, which represents the
	response of algae to light, and mathematically, it is deﬁned as
	
	\begin{equation*}
		M(G)=\left\{\begin{array}{ll}\geq 0, & \mbox{when }~~ G\leq G_{c},\\< 0, & \mbox{when }~~G>G_{c}.  \end{array}\right. 
	\end{equation*}
	
	\section{Governing equations}
	
	In common with the previous models of bioconvection, we
	assume a monodisperse cell population that can be modelled by a continuous distribution. The suspension is diluted so that the volume fraction of the cells is small and cell–cell interactions are negligible. Each cell has a volume $\vartheta$ and density $\rho + \Delta\rho$, $\rho$ is the density of the ﬂuid in which the cells swim and $\Delta\rho\ll \rho$. Let $\boldsymbol{u}$ and $n$, respectively, denote the ﬂuid velocity and concentration in the suspension. Supposing that the suspension is incompressible, conservation of mass implies
	
	\begin{equation}\label{6}
	\boldsymbol{\nabla}\cdot \boldsymbol{u}=0.
	\end{equation}
	We assume that the Stokeslet due to the negative buoyancy of
	the cells is dominant and all other contributions of the cells to the bulk stress are negligible. Thus, neglecting all effects on the ﬂuid except the cells’ negative buoyancy, the momentum equation under the
	Boussinesq approximation is
	
	\begin{equation}\label{7}
	\rho\left(\frac{D\boldsymbol{u}}{D t}\right)=-\boldsymbol{\nabla} P_e+\mu{\nabla}^2\boldsymbol{u}-n\vartheta g\Delta\rho\hat{\boldsymbol{z}}-\rho g(1-\beta(T-T_0))\hat{\boldsymbol{z}}.
	\end{equation}
	
	where $P_e$ is the excess pressure above the hydrostatic pressure, $g$ is the acceleration due to gravity, and $\mu$ is the viscosity of the suspension, which is assumed to be that of the ﬂuid. Here, $D/Dt$ is the material time derivative and $\beta$ is the coefficient of thermal expansion.
	
	The equation for cell conservation is
	\begin{equation}\label{8}
	\frac{\partial n}{\partial t}=-\boldsymbol{\nabla}\cdot \boldsymbol{F},
	\end{equation}
	
	where the ﬂux of cells is
	\begin{equation*}
		\boldsymbol{F}=n(\boldsymbol{u}+W_c<\boldsymbol{p}>)-\boldsymbol{D}\cdot\boldsymbol{\nabla} n.
	\end{equation*}
	
	The ﬁrst term on the right-hand side of Eq. (7) is the ﬂux due to the advection of the cells by the bulk ﬂuid ﬂow. The second term arises due to the average swimming velocity of the cells.
	The third term represents the random component of the cell locomotion. We choose the diffusivity tensor $\boldsymbol{D}$ to be isotropic and constant. Thus, $\boldsymbol{D}=DI$, where $D$ is the diffusion coefﬁcient and $I$ is the identity tensor. Thus, the ﬂux of cells becomes
	
	\begin{equation*}
		\boldsymbol{F}=n(\boldsymbol{u}+W_c<\boldsymbol{p}>)-{D}\cdot\boldsymbol{\nabla} n.
	\end{equation*}
	
	Here, each algal cell is idealized as a homogeneous spherical body having a uniform distribution of mass and purely phototactic, and thus, it is centre of gravity and centre of buoyancy coincide (i.e., the torque due to gravity does not exist). Two major assumptions have been made in deriving the cell ﬂux vector in the proposed model. First, since the
	algal cells are purely phototactic, the effect of viscous torque due to shear in the ﬂow, which might contribute to the horizontal component of the mean swimming orientation, is neglected. It can be justiﬁed clearly from the linear stability theory that the perturbed vorticity, and hence the viscous torque, is inﬁnitesimally small. Hence, the viscous reorientation can be neglected on the limit. Second, the diffusion tensor, which should be a function of light intensity, is assumed to be a constant isotropic tensor, instead of deriving it from a swimming velocity autocorrelation function. Thus, this model can be considered
	to be valid in the limiting case, to understand the complexities of bioconvection due to phototaxis before exploring more detailed complex models. These assumptions allow us to eliminate the Fokker–Planck equation (which relates the degree of alignment of phototactic algae with the ﬂuid below) from the governing system for phototactic bioconvection.
	
	Also, the thermal energy equation
	\begin{equation}\label{9}
	\rho c\big[\frac{\partial T}{\partial t}+\boldsymbol{\nabla}\cdot(\boldsymbol{u}T)\big] =\alpha{\boldsymbol{\nabla}}^2 T,
	\end{equation}	
	where $\rho c$ is the volumetric heat capacity of water, and $\alpha$ is the thermal conductivity of water.
	
	\section{Boundary conditions}
	The lower boundary ($z=0$) is taken as rigid, while the upper
	boundary ($z=H$) may be stress-free or rigid. Since there is no ﬂux of cells at each boundary, the boundary conditions are
	
	\begin{subequations}
		\begin{equation}\label{10a}
		\boldsymbol{u}=\boldsymbol{F}\cdot\hat{\boldsymbol{z}}=0\qquad \text{on}~~~~ z=0,
		\end{equation}
		\begin{equation}\label{10b}
		\boldsymbol{u}\cdot\hat{\boldsymbol{z}}=\boldsymbol{F}\cdot\hat{\boldsymbol{z}}=0\qquad \text{on}~~~~ z=H,
		\end{equation}
		for rigid boundaries
		\begin{equation}\label{10c}
		\boldsymbol{u}\times\hat{\boldsymbol{z}}^*=0\qquad \text{on}~~~~ z=0,H,
		\end{equation}
		while for a free boundary
		\begin{equation}\label{10d}
		\frac{\partial^2}{\partial z^2}(\boldsymbol{u}\times\hat{\boldsymbol{z}})=0\qquad \text{on}~~~~ z=H,
		\end{equation}	
		for temperature on both boundaries
		\begin{equation}\label{10e}
		T=T_0+\Delta T\qquad \text{on} ~~~~z=0,
		\end{equation}
		\begin{equation}\label{10f}
		T=T_0\qquad \text{on}~~~~ z=H.
		\end{equation}
	\end{subequations}
	
	The top boundary is exposed to a uniform oblique solar irradiation of magnitude $I_t$. The top and the bottom boundaries are also assumed to be non-reﬂecting, thus,
	
	\begin{equation*}
		\text{at}~z=H,~~~~~	I(x,y,z,\theta,\phi)=I_t\delta(\boldsymbol{s},\boldsymbol{s}_0)~~~~~~~~(\pi/2\leq\theta\leq\pi).
	\end{equation*}
	\begin{equation*}
		\text{at}~z=0,~~~~~I(x,y,z,\theta,\phi)=0~~~~~~~~(0\leq\theta\leq\pi/2).
	\end{equation*}
	
	\section{Scaling of the equations}
	
	The governing system for bioconvection can be made dimensionless by scaling all lengths on $H$, the depth of the layer, time on diffusive timescale $H^2/\alpha_f$, and the bulk ﬂuid velocity on $\alpha_f/H$. The appropriate scaling for the pressure is $\mu \alpha_f/H^2$; the cell concentration is scaled $\bar{n}$, the mean concentration; and temperature is scalled by $(T-T_0)/\Delta T$. For convenience, we keep the same on $n$ notation for the dimensional and nondimensional variables. The recasting of the governing system for bioconvection in terms of the nondimensional variables is as follows
	
	\begin{equation}\label{11}
	\boldsymbol{\nabla}\cdot\boldsymbol{u}=0,
	\end{equation}
	\begin{equation}\label{12}
	P_r^{-1}\left(\frac{D\boldsymbol{u}}{D t}\right)=-\nabla P_{e}+\nabla^{2}\boldsymbol{u}-R_bn\hat{\boldsymbol{z}}-R_m\hat{\boldsymbol{z}}+R_TT\hat{\boldsymbol{z}},
	\end{equation}
	\begin{equation}\label{13}
	\frac{\partial{n}}{\partial{t}}=-\boldsymbol{\nabla}\cdot\boldsymbol{F},
	\end{equation}
	where the ﬂux of cells is
	\begin{equation}\label{14}
	\boldsymbol{F}=\boldsymbol{n{\boldsymbol{u}}+\frac{1}{Le}nV_{c}<{\boldsymbol{p}}>-\frac{1}{Le}{\boldsymbol{\nabla}}n},
	\end{equation}
	\begin{equation}\label{15}
	\frac{\partial T}{\partial t}+\boldsymbol{\nabla}\cdot(\boldsymbol{u}T) =\boldsymbol{\nabla}^2 T.
	\end{equation}
	
	Here, $P_r=\nu/\alpha_f$ is the Prandl number, $V_c=W_sH/D$ is the scaled swimming speed, $R_b=\bar{n}v g\Delta{\rho}H^{3}/\mu \alpha_f$ is the bioconvective Rayleigh number, $R_T= g\beta\Delta{T}H^{3}/\mu\alpha_f$ is the thermal Rayleigh number and $R_m=\rho gH^3/\mu \alpha_f$ is the basic density Rayleigh number and Lewis number $Le=\alpha_f/D$. Here, $\alpha_f$ is the thermal diffusivity of water. In dimensionless form, the boundary conditions become
	
	\begin{subequations}
		\begin{equation}\label{16a}
		\boldsymbol{u}=\boldsymbol{F}\cdot\hat{\boldsymbol{z}}=0\qquad \text{on}~~~~ z=0,
		\end{equation}
		\begin{equation}\label{16b}
		\boldsymbol{u}\cdot\hat{\boldsymbol{z}}=\boldsymbol{F}\cdot\hat{\boldsymbol{z}}=0\qquad \text{on}~~~~ z=1,
		\end{equation}
		for rigid boundaries
		\begin{equation}\label{16c}
		\boldsymbol{u}\times\hat{\boldsymbol{z}}=0\qquad \text{on}~~~~ z=0,1,
		\end{equation}
		while for a free boundary
		\begin{equation}\label{16d}
		\frac{\partial^2}{\partial z^2}(\boldsymbol{u}\times\hat{\boldsymbol{z}})=0\qquad \text{on}~~~~ z=1,
		\end{equation}	
		for temperature on both boundaries
		\begin{equation}\label{16e}
		T=T_0+\Delta T\qquad \text{on} ~~~~z=0,
		\end{equation}
		\begin{equation}\label{16f}
		T=T_0\qquad \text{on}~~~~ z=1.
		\end{equation}
	\end{subequations}
	
	In terms of the nondimensional variables, the RTE [see Eq. (2)]
	becomes
	\begin{equation}\label{17}
	\boldsymbol{s}\cdot\nabla I(\boldsymbol{x},\boldsymbol{s})+\varkappa nI(\boldsymbol{x},\boldsymbol{s})=0
	\end{equation}
	where $\kappa_H=\kappa\bar{n}H$ is the (vertical) optical depth of the suspension. In dimensionless form, the intensity at the top and bottom becomes

	\begin{equation}\label{18}
	I(x, y, z = 1, \theta, \phi) = I_t\delta(\boldsymbol{s}-\boldsymbol{s_0}),\qquad (\pi/2\leq\theta\leq\pi),
	\end{equation}
	\begin{equation}\label{19}
	I(x, y, z = 0, \theta, \phi) =0,\qquad (0\leq\theta\leq\pi/2),
	\end{equation}
	
	\section{The basic solution}
	Equations (\ref{11})–(\ref{15}) and Eq. (\ref{17}) together with the boundary
	conditions possess a static equilibrium solution in which
	
	\begin{equation*}
		\boldsymbol{u}=0,~n=n_s(z),~T=T_s(z)~\text{and}~ I=I_s(z,\theta).
	\end{equation*}
	
	The total intensity $G_s$ at the basic state is given by  
	
	\begin{equation*}
		G_s=\int_0^{4\pi}I_s(z,\theta)d\Omega.
	\end{equation*}
	
	If horizontal homogeneity and azimuthally isotropic are assumed
	(as considered here), the RTE becomes
	
	\begin{equation}\label{20}
	\frac{dI_s(z,\theta)}{dz}+\frac{\tau_H n_sI_s(z,\theta)}{\cos\theta}=0,
	\end{equation}
	
	subject to the top boundary condition
	
	\begin{equation}\label{21}
	I_s(1,\theta) =I_t\delta(\boldsymbol{s}-\boldsymbol{ s}_0).
	\end{equation}
	
	Solving Eqs. (\ref{20}) and (\ref{21}), we get
	
	\begin{equation}\label{22}
	dI_s(z,\theta)=I_t\delta(\boldsymbol{s}-\boldsymbol{ s}_0)\exp\bigg(\int_z^1\frac{\tau_Hn_s(z')}{\cos\theta}dz'    \bigg),
	\end{equation}
	
	The basic total intensity is written as
	
	\begin{equation}\label{23}
	G_s(z)=\int_0^{4\pi}I_s(z,\theta)d\Omega=I_t\exp\bigg(\frac{-\int_z^1\tau_Hn_s(z')}{\cos\theta_0}dz'\bigg).
	\end{equation}
	
	Now, the radiative heat ﬂux, $\boldsymbol{q}_s(z)$, in basic state is given by
	
	\begin{equation*}
		\boldsymbol{q}_s(z)=\int_0^{4\pi}I_s(z,\theta)\hat{\boldsymbol{s}}d\Omega.
	\end{equation*}
	The radiative heat ﬂux normal to the boundary surfaces of the algal suspension is given by
	
	\begin{equation*}
		(\boldsymbol{q_s}\cdot\hat{\boldsymbol{z}})\hat{\boldsymbol{z}}=\bigg[I_t\exp\left(\frac{-\tau_H\int_z^1n_s(z')dz'}{\cos\theta_0}\right)\cos\theta_0\bigg](-\hat{\boldsymbol{z}})
	\end{equation*}	
	
	Since the illuminating source lies in the opposite direction to the radiative heat ﬂux vector, the mean swimming direction in the basic state is expressed as	
	
	\begin{equation*}
		<\boldsymbol{p}_s>=M_s\hat{\boldsymbol{z}}~~\text{where}~~M_s=M({G}).
	\end{equation*}
	
	The basic concentration, n s ðzÞ, satisﬁes
	
	\begin{equation}\label{24}
	\frac{dn_s}{dz}-V_cM_sn_s=0,
	\end{equation}
	
	which is supplemented by the cell conservation relation
	
	\begin{equation}\label{25}
	\int_0^1n_s(z)dz=1.
	\end{equation}
	
	The basic temperature $T_s(z)$ is satisfied satisfy
	
	\begin{equation}\label{26}
	\frac{d^2T_s}{dz^2}=0.
	\end{equation}
	where 
	
	\begin{equation}\label{27}
	T_s-1=0,~~\text{at}~~z=0
	\end{equation}
	\begin{equation}\label{28}
	T_s=0,~~\text{at}~~z=1.
	\end{equation}
	
	By introducing a new variable $\varpi=\int_1^zn_sdz$, Eq.~(\ref{25}) becomes
	\begin{equation}\label{29}
	\frac{d^2\varpi}{dz^2}-V_cT_s\frac{d\varpi}{sz}=0,
	\end{equation}
	with boundary condition
	
	\begin{equation}\label{30}
	\varpi+1=0,~~\text{at}~~z=0
	\end{equation}
	\begin{equation}\label{31}
	\varpi=0,~~\text{at}~~z=1.
	\end{equation}
	
	In terms of $\varpi$, the basic total intensity $G_s$ can be expressed as
	\begin{equation*}
		G_s(z)=I_t\exp\bigg(\frac{\tau_H\varpi}{\cos\theta_0}\bigg)
	\end{equation*}
	
	Equations (\ref{26}) and (\ref{29}) with specific boundary conditions constitute a boundary value problem, which is solved numerically using a shooting method.
	
	\section{Linear stability of the problem}
	
	To perform linear stability analysis, the basic state is perturbed via inﬁnitesimal disturbances as
	
	\begin{equation}\label{32}
	\begin{pmatrix}
	\boldsymbol{u}\\n\\T\\<\boldsymbol{I}>
	\end{pmatrix}
	=
	\begin{pmatrix}
	0\\n_s\\T_s\\<\boldsymbol{I}_s>
	\end{pmatrix}
	+\epsilon
	\begin{pmatrix}
	\boldsymbol{u}_1\\n_1\\T_1\\<\boldsymbol{I}_1>
	\end{pmatrix}
	+O(\epsilon^2),
	\end{equation}
	
	where $\boldsymbol{u}_1=(u_1,v_1,w_1)$.
	
	The linearized equations about the basic state are
	
	\begin{equation}\label{33}
	\boldsymbol{\nabla}\cdot \boldsymbol{u}_1=0,
	\end{equation}
	
	\begin{equation}\label{34}
	P_r^{-1}\left(\frac{\partial \boldsymbol{u_1}}{\partial t}\right)=-\boldsymbol{\nabla} P_{e}+\nabla^{2}\boldsymbol{u}_1-R_bn_1\hat{\boldsymbol{z}}+R_TT_1\hat{\boldsymbol{z}},
	\end{equation}
	
	\begin{equation}\label{35}
	\frac{\partial{n_1}}{\partial{t}}+\frac{1}{Le}V_c\boldsymbol{\nabla}\cdot(<\boldsymbol{p_s}>n_1+<\boldsymbol{p_1}>n_s)+w_1\frac{dn_s}{dz}=\frac{1}{Le}\boldsymbol{\nabla}^2n_1,
	\end{equation}
	
	\begin{equation}\label{36}
	\frac{\partial{T_1}}{\partial t}-w_1\frac{dT_s}{dz}=\boldsymbol{\nabla}^2T_1.
	\end{equation}
	
	Now the total intensity, $G$, is given by
	
	\begin{equation*}
		G_s(z)=I_t\bigg[\frac{\int_1^z(\tau_Hn_s(z)+\epsilon n_1 +O(\epsilon^2))dz'}{\cos\theta_0}\bigg]
	\end{equation*}
	
	The steady intensity is perturbed, and after simpliﬁcation, we get
	
	\begin{equation*}
		G=G_s+\epsilon G_1+O(\epsilon^2)
	\end{equation*}
	
	where
	
	\begin{equation}\label{37}
	G_1=I_t\exp\left(\frac{\int_1^z\tau_H n_s(z')dz'}{\cos\theta_0}\right)\left(\frac{\int_1^z\tau_H n_1 dz'}{\cos\theta_0}\right),
	\end{equation}
	
	Hence, the perturbed mean swimming orientation [i.e., $T(G)\hat{\boldsymbol{z}}$] at $O(\epsilon^2)$ for a non-scattering algal suspension is expressed as
	
	\begin{equation}\label{38}
	<\boldsymbol{p_1}>=G_1\frac{dM_s}{dG}\hat{\boldsymbol{z}}.
	\end{equation}
	
	We eliminate $P_e$ and the horizontal component of $\boldsymbol{u_1}$ by taking the curl of Eq. (\ref{34}) twice and retaining the z-component of the result. Then, Eqs. (\ref{33})–(\ref{37}) reduce to two equations for $w_1$ and $n_1$. Now, this equation can be decomposed into normal mode as
	
	\begin{equation}\label{39}
	\begin{pmatrix}
	w_1\\n_1\\T_1
	\end{pmatrix}
	=
	\begin{pmatrix}
	W(z)\\\Theta(z)\\T(z)
	\end{pmatrix}
	+\exp{[\sigma t+i(k_xx+k_xy)]},
	\end{equation} 
	
	$W(z)$, $\Theta(z)$, and $T(z)$ represent the variations in the $z$ direction, while $k_x$ and $k_y$ are the horizontal wavenumbers. The complex growth rate of the disturbances is denoted by $\sigma$.  
	
	The linear stability equations become
	\begin{equation}\label{40}
	\left(\sigma P_r^{-1}+k^2-\frac{d^2}{dz^2}\right)\left( \frac{d^2}{dz^2}-k^2\right)W=R_bk^2\Theta-R_Tk^2T,
	\end{equation}
	\begin{equation}\label{41}
	\left(\sigma Le+k^2-\frac{d^2}{dz^2}\right)\Theta+\aleph_0(z)\int_z^1 \Theta dz+\aleph_1(z)\Theta+\aleph_2(z)\frac{d\Theta}{dz}=-Le\frac{dn_s}{dz}W,
	\end{equation}
	
	\begin{equation}\label{42}
	\left(\frac{d^2}{dz^2}-k^2-\gamma\right)T(z)=\frac{dT_s}{dz}W(z),
	\end{equation} 
	where 
	\begin{subequations}
		\begin{equation}\label{43a}
		\aleph_0(z)=-(\tau_H/\cos\theta_0) V_c\frac{d}{dz}\left(n_sG_s^c\frac{dM_s}{dG}\right),
		\end{equation}
		\begin{equation}\label{43b}
		\aleph_1(z)=2(\tau_H/\cos\theta_0) V_c n_s G_s\frac{dM_s}{dG},
		\end{equation}
		\begin{equation}\label{43c}
		\aleph_2(z)=V_cM_s.
		\end{equation}
	\end{subequations}  
	
	subject to the boundary conditions
	\begin{equation}\label{44}
	W=\frac{dW}{dz}=\frac{d\Theta}{dz}-\aleph_2(z)\Theta-n_sV_c(\tau_H/\cos\theta_0)G_s\bigg(\int_z^1\Theta dz\bigg)\frac{dM_s}{dG}=0,~~\text{at}~~~z=0,
	\end{equation}
	\begin{equation}\label{45}
	W=\frac{d^2W}{dz^2}=\frac{d\Theta}{dz}-\aleph_2(z)\Theta=0,~~\text{at}~~~z=1.
	\end{equation}
	At a rigid upper surface, the condition in Eq. (\ref{45}) is replaced by
	\begin{equation}\label{46}
	W=\frac{dW}{dz}=\frac{d\Theta}{dz}-\aleph_2(z)\Theta=0,~~\text{at}~~~z=1.
	\end{equation}
	Here, $k=\sqrt{k_x^2+k_y^2}$ is the horizontal wavenumber and it represents the modulation of the bioconvection pattern in horizontal directions. Equations (\ref{40})-(\ref{42}) form an eigenvalue problem for $\sigma$ as a function of the dimensionless parameters $\theta_0$, $k$, $P_r$, $V_c$, $\tau_H$, $R_b$ and $R_T$. The basic
	state becomes unstable whenever $Re(\sigma)>0$.
	
	Introducing a new variable
	
	\begin{equation}\label{47}
	\Phi(z)=\int_z^1\Theta(z')dz',
	\end{equation} 
	
	the linear stability equations become (writing $D= d/dz$)
	
	\begin{equation}\label{48}
	D^4W-(2k^2+\sigma Le P_r^{-1})D^2W+k^2(k^2+\sigma Le P_r^{-1})W=-R_bk^2D\Phi+R_Tk^2T(z),
	\end{equation}
	\begin{equation}\label{49}
	D^3\Phi-\aleph_3(z)D^2\Phi-(\sigma Le+k^2+\aleph_2(z))D\Phi-\aleph_1(z)\Phi-\aleph_0(z)=Le Dn_sW, 
	\end{equation}
	\begin{equation}\label{50}
	\left(D^2-k^2-\sigma\right)T(z)=DT_sW,
	\end{equation}
	with
	
	\begin{equation}\label{51}
	W=DW=D^2\Phi-\aleph_2(z)D\Phi-n_sV_c(\tau_H/\cos\theta_0)G_s\frac{dM_s}{dG}\Phi=0,~~\text{at}~~~z=0,
	\end{equation}
	\begin{equation}\label{52}
	W=D^2W=D^2\Phi-\aleph_2(z)D\Phi=0,~~\text{at}~~~z=1.
	\end{equation}
	At a rigid upper surface, the condition in Eq. (37) is replaced by
	\begin{equation}\label{53}
	W=DW=D^2\Phi-\aleph_2(z)D\Phi=0,~~\text{at}~~~z=1.
	\end{equation}
	\begin{equation}\label{54}
	T(z)=0,~~\text{at}~~z=0,1.
	\end{equation}
	
	There is also an extra boundary condition
	
	\begin{equation}\label{55}
	\Phi(z)=0,~~~ \text{at}~~~z=1.
	\end{equation}
	
	which follows from Eq. (\ref{47}). The new aspects of the proposed model via the effects of oblique collimated irradiation are incorporated in the terms $\aleph_2$, $\aleph_2$, and $\aleph_2$, which are found on the perturbed cell conservation equation, Eq. (\ref{49}). 
	
	
	\section{SOLUTION PROCEDURE}
	
	Numerical solutions to Eqs. (\ref{48}) - (\ref{50})) with appropriate boundary conditions are calculated with a fourth-order accurate, ﬁnite-difference scheme based on Newton–Raphson–Kantorovich (NRK) iterations. This NRK routine solves a system of ﬁrst-order, nonlinear, coupled, ordinary differential equations (ODEs) resulting from the two-point boundary value problem. In the NRK routine, a guess is required for the solution of the system of ordinary differential equations (ODEs). Next, the guess is reﬁned by successive iterations until the desired accuracy is reached. The numerical results are checked for accuracy by comparing the solution computed by the NRK routine (used in this study) of a benchmark problem in the same vicinity with another (solution) calculated by a routine, which is different from the former (NRK routine), such as the solution of the
	mixed-convection ﬂow in a vertical pipe. It reveals that the solution generated by the said routine (NRK) is in excellent agreement with the one, which utilizes the Chebyshev collocation method (see Table I of Ref. 39 for details). Furthermore, the numerical scheme is tested for many different parameters and mesh sizes. It shows that the solutions obtained using the scheme for the same parameters on different meshes always agree to 5 (or more) signiﬁcant ﬁgures for a minimum of 51 mesh points. This numerical scheme (NRK routine) is utilized to
	calculate the neutral stability curves (or neutral curves) in the $(k, R)$-plane, where $R=(R_b, R_T)$ or the growth rate, $Re(\sigma$), as a function of $R$ for a ﬁxed set of other parameters. Initially, the values of $P_r$, $V_c$, $\tau_H$, $\theta_0$, and $k$ are supplied, and the values of $W$, $\Phi$ and $T$ are estimated either from the previous numerical results or by imposing sinusoidal variation in $W$, $\Phi$ and $T$. Once a solution is obtained, this solution can be used as an initial guess for the neighbouring parameter values.
	
	For a given set of other parameters, there are an inﬁnite number of branches of the neutral curve $R^{(n)}(k)$, with $n=1,2,3...$; each one representing a particular solution of the linear stability problem. The solution branch of most interest is the one on which $R$ has its minimum value $R^c$ and the corresponding bioconvective solution, that is, $(k^c, R^c)$, where $R^c$ is either $R_b^c$ or $R_T^c$ is called the most unstable solution. The wavelength of the initial disturbance can be calculated from it by using the relation $\lambda_c=2\pi/k^c$. Please note here that the linear stability theory is expected to explain the initial bioconvection pattern spacing before the nonlinear effects in the system become dominant. Solutions consist of convection cells stacked one above another along the depth of the suspension. A solution is said to be of mode $n$ if it has n convection cells stacked vertically one on
	another. In many instances, the most unstable solution occurs on the $R^{(n)}(k)$ branch of the neutral curve and it is mode 1. 
	
	A neutral curve is deﬁned as the locus of points where $Re(\sigma)=0$. If in addition $Im(\sigma)=0$ on such a curve, then the principle of exchange of stabilities is said to be valid and the bioconvective solution is called stationary (non-oscillatory). Alternatively, if $Im(\sigma)\neq 0$ then oscillatory solutions exist. If the most unstable solution remains on the oscillatory branch of the neutral curve, then the solution is called overstable. When there is competition between the stabilizing
	and destabilizing processes, the oscillatory solution evolves. In this instance, a single oscillatory branch bifurcates from the stationary branch of the corresponding neutral curve at some $k=k_b$ and exists for $k<k_b$.
	
	To estimate the parameters required for the present study, we assume that we are dealing with a phototactic microorganism similar to Chlamydomonas. To retain the resulting model to be more rational and can be comparable with earlier studies on (phototactic) bioconvection, we use the same parameter values as taken by Refs. 17, 22, and 40. Accounting for refraction at the air-water interface, the range of the declination angle $\theta_0$ is restricted such that $0.661\leq\cos\theta_0\leq 1.0$; and this implies that $0\leq\theta_i(deg)\leq 48.6$ approximately. Following Daniel et al., we have estimated the approximate range of the angle of incidence $\theta_i$ as $0\leq\theta_i\leq 80$ for our proposed model. The radiation characteristics required here are calculated as given in Ghorai and Panda. Thus, for a 0.50cm deep suspension, optical depth $\tau_H$ varies in the range from 0.25 to 1. The scaled swimming speed can be calculated for a 0.5 and 1.0cm deep suspension as $V_c=10$ and $V_c=20$, respectively.
	
	\section{NUMERICAL RESULTS}
	We have systematically investigated the effects of the angle of
	incidence $\theta_i$ on the base concentration proﬁles and the corresponding neutral curves by varying it via discrete values, keeping the other governing parameters ﬁxed. The values of $Pr=5$ and $I_t=0.8$ are kept ﬁxed throughout to compare our model with the existing rational models on (phototactic) bioconvection. The representative values of the cell swimming speed and the
	optical (vertical) depth are $V_c=10$, 15, 20, and $\tau_H=0.5$, 1.0, respectively. The value of the critical total intensity $G_c$ is selected such that
	the maximum basic concentration occurs at around mid-height of the suspension domain, while the angle of incidence remains ﬁxed as $\theta_i=0$ (i.e., vertical collimated irradiation case). Then, the angle of incidence $\theta_i\in[0:80]$ is varied discretely with an increment of $10$ starting from $0$ and we study their effects on the bioconvection patterns at instability. Please note that due to a large number of parameter values, it is very difficult to obtain a comprehensive picture across the whole parameter domain. Thus, we have taken a discrete set of ﬁxed parameter values to study their effects on the onset of bioconvection. We determine the bioconvective Rayleigh number, $R_b$, at the onset of bioconvection as a function of wavenumber, $k$, for various thermal Rayleigh numbers, $R_T$.
	We also divide the obtained results into two categories based on the suspension top boundary condition (e.g., stress-free or rigid). Then, we address them separately by taking two cases.
	
	\section{Case I: Stress-free upper surface}
	To study the effects of angle of incidence on the onset of bioconvection, here we consider the case when the upper surface of the algal
	suspension is stress-free.

	\begin{figure}[!htbp]
		\centering
		\includegraphics[width=10cm]{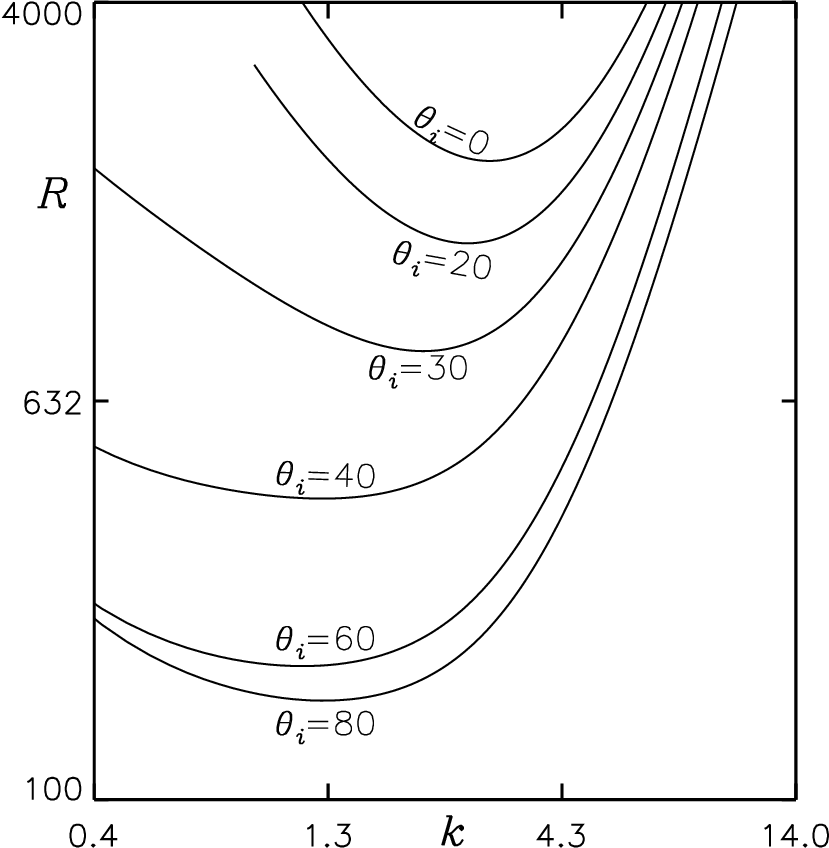}
		\caption{\footnotesize{The Neutral curves for various values of the angle of incidence $\theta_i$. Here, the upper surface is stress-free and the other governing parameters $V_c=10$, $\tau_H=0.5$ and $R_T=50$ are fixed.}}
		\label{fig2}
	\end{figure}
	In Fig.~\ref{fig2}, we show the neutral curves at the variation in the angle of incidence as $\theta_i=0$, 20, 30, 40, 60, and 80, respectively, while
	keeping the other governing parameters $V_c=10$; $\tau_H=0.5,$, $R_T=50$ and $G_c=0.63$ ﬁxed. When $\theta_i=0$, the maximum basic concentration occurs
	at around mid-height of the domain resulting into a horizontal concentrated sublayer. The region below (above) the sublayer is gravitationally unstable (locally stable). As the value of $\theta_i$ is increased up to 80, the location of the maximum basic concentration (i.e., sublayer) shifts toward the top of the domain, and simultaneously, the value of maximum basic concentration increases. The width of the upper (lower) stable (unstable) region gradually decreases (increases) at the same time. Thereby, the buoyancy of the upper stable ﬂuid does not inhibit to convection adequately, whereas the lower unstable convective region supports it very strongly due to a substantial increase in its height. Eventually, both the critical wavenumber $k^c$ and bioconvective Rayleigh number $R_b^c$ at bioconvective instability decrease and thus, the suspension becomes more unstable as the angle of incidence $\theta_i$ is increased to higher non-zero values [see Fig. \ref{fig2}].
	
	\begin{figure}[!htbp]
		\centering
		\includegraphics[width=10cm]{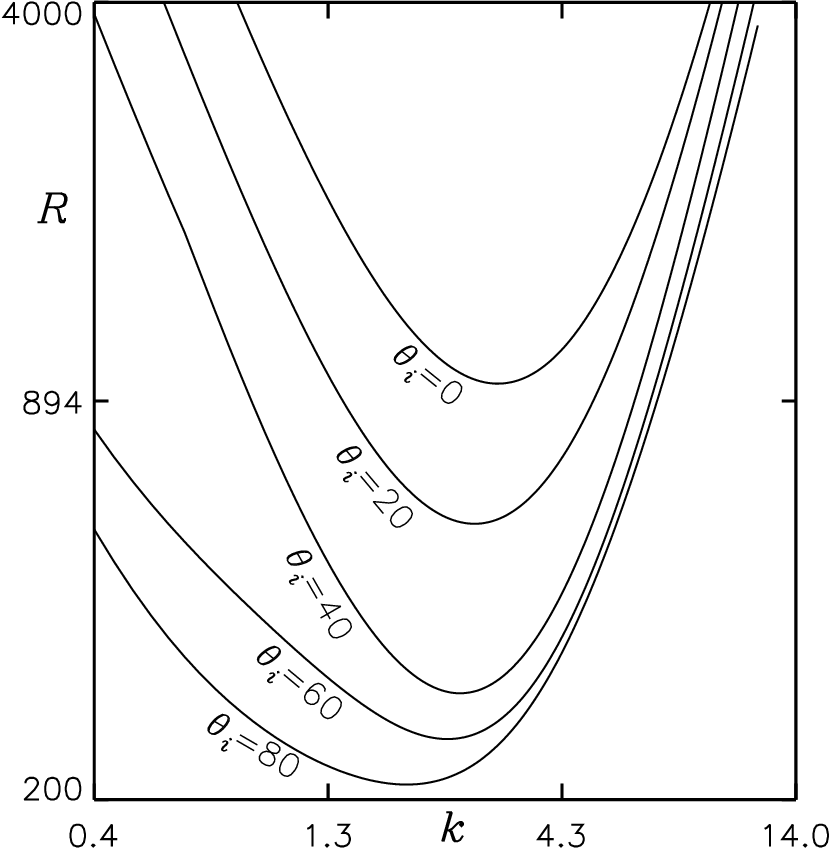}
		\caption{\footnotesize{The Neutral curves for various values of the angle of incidence $\theta_i$. Here, the upper surface is rigid and the other governing parameters $V_c=10$, $\tau_H=1$ and $R_T=50$ are fixed.}}
		\label{fig3}
	\end{figure}
	
	\section{Case-II: Rigid upper surface}
	To emphasize the comparison between the theoretical predictions of the proposed model on bioconvection and the corresponding quantitative studies, we analyze here the solutions at the instability in a suspension of phototactic algae contained between two rigid surfaces located at $z=0$ and $z=1$, respectively. It may also be interesting to know that the collection of cells at the top surface may form a rigid-like packed layer if the upper boundary is open to the air medium.

	Figure \ref{fig3} shows the neutral curves at the variation in the angle
	of incidence as $\theta_i=0$, 20, 40, 60, and 80, respectively, while the
	other governing parameters $V_c=10$; $\tau_H=01,$, $R_T=50$ and $G_c=0.495$ ﬁxed. remain constant. At $\theta_i=0$, the location of the maximum basic concentration is around the mid-height of the chamber and the corresponding stationary branch of the neutral curve possesses the most unstable solution leading the bioconvective solution to be non-oscillatory (stationary). As $\theta_i$ is increased further to higher non-zero values, the location of the maximum basic concentration shifts toward the top of the chamber and simultaneously the value of maximum concentration increases. The most unstable solution tends to remain in the stationary branch of the neutral curve as $\theta_i$ is increased up to 80. The width of the upper (lower) stable (unstable) region gradually decreases (increases) at the same time. Thereby, the buoyancy of the upper stable ﬂuid does not inhibit convection adequately, whereas the lower unstable convective region supports it very strongly due to a substantial increase in its height. Eventually, both the critical wavenumber $k^c$ and bioconvective Rayleigh number $R_b^c$ at bioconvective instability decrease and thus, the suspension becomes more unstable as the angle of incidence $\theta_i$ is increased to higher non-zero values [see Fig. \ref{fig3}].
	
	\section{Conclusion}
	
	In this innovative model of thermal bioconvection driven by thermal effect and phototaxis, we investigate the emergence of bio-thermal convection within a suspension comprising non-scattering phototactic algae. The oblique collimated irradiation strikes the suspension top at a ﬁxed angle of incidence. The linear stability of the suspension has been analyzed using this model. In the basic state, it is observed that self-shading becomes dominant (via increment in slant-path length) at the increment in angle of incidence, and thus, the algae receive light of low intensities (via variation in angle of incidence) at a ﬁxed interior depth. Thereby, the location of maximum basic concentration shifts toward the suspension top and the value of maximum basic concentration increases at the increment in angle of incidence. When the value of optical depth is higher, the base concentration proﬁles become steeper as the angle of incidence is increased, while keeping the other governing parameters ﬁxed.
	
	Linear stability analysis predicts that the perturbation to the basic
	steady state is either stationary or oscillatory (overstable). Furthermore, the critical Rayleigh number usually decreases as the angle of incidence
	is increased. As a result, the algal suspension becomes more unstable
	at the increment in angle of incidence.

	
	\section*{ Availability of Data}
	The supporting data of this article is available within the article. 
	\nocite{*}
	\section*{REFERENCES}
	\bibliography{Oblique}
	
\end{document}